\documentclass[11pt,reqno]{amsart}
\usepackage{amssymb,latexsym}
\usepackage{graphicx}
\usepackage{fancyhdr,amssymb}
\usepackage{url}		

\theoremstyle{plain}
\numberwithin{equation}{section}

\begin{document}
\fancyhead{}
\renewcommand{\headrulewidth}{0pt}
\fancyfoot{}
\fancyfoot[LE,RO]{\medskip \thepage}

\setcounter{page}{1}

\title[Fibonacci identities and  Fibonacci pairs
 ]{Fibonacci identities and  Fibonacci pairs}
\author{Cheng Lien Lang}
\address{Department of Applied  Mathematics\\
                I-Shou University\\
                Kaohsiung, Taiwan\\
                }
\email{cllang@isu.edu.tw}
\thanks{}
\author{Mong Lung Lang}
\address{Singapore
669608,
Singapore   }
\email{lang2to46@gmail.com}

\begin{abstract}
 A Fibonacci pair   $F_s(w,x)$ of rank $s$ is a pair of $s\times s$ nonsingular matrices such that
  $wx=xw$ and that  entries of $aw^n$ and $axw^m$
   are polynomials of Fibonacci or Lucas numbers  for some $a\ne 0$.
We construct identities
 systematically by the study of $F_2(w,x)$ and $F_3(w, x)$.

 \end{abstract}

\maketitle
\vspace{-.8cm}

\section{Introduction}
A Fibonacci pair   $F_s(w,x)$ of rank $s$ is a pair of $s\times s$ nonsingular matrices such that
  $wx=xw$ and that  entries of $aw^n$ and $axw^m$
   are polynomials of Fibonacci or Lucas numbers
    for some $a\ne 0$.
Existence of $F_s(w, x)$  such that $\left < \bar w, \bar
x\right > \subseteq GL(s, \Bbb C)/\left < cI_s\,:\, c\ne 0\right > $  is not cyclic is
   guaranteed by the following.

   \smallskip
   \noindent {\bf Theorem 1.1.} {\em  Let
   $w={\tiny \left (\begin{array}{cc}
1&1\\
1&0\\
\end{array}\right)}$ and $x ={\tiny \left (\begin{array}{cc}
1&2\\
2&-1\\
\end{array}\right)}$. Then $wx=xw$,
$\left < \bar w, \bar
x\right > $ is not cyclic,
$$w^n= \left (\begin{array}{cc}
F_{n+1}&F_n\\
F_n&F_{n-1}\\
\end{array}\right) \mbox{  and }
xw^n= \left (\begin{array}{cc}
L_{n+1}&L_n\\
L_n&L_{n-1}\\
\end{array}\right).\eqno(1.1)$$}
Fibonacci pairs can be constructed easily (see Proposition 2.1 and 5.1).
 In this article,
we give 2  Fibonacci pairs  of rank 2 (subsection 1.1)
 and 2 pairs of rank 3 (Section 5).
  Note that these pairs are  constructed in such a way that
     the groups $\left < \bar w,\bar  x\right >\subseteq
   GL(s, \Bbb C)/\left < cI_s\,:\, c\ne 0\right >
  $ are not cyclic.
For each  pair of rank 2, we  construct
   two types of identities,
   matrix identities (subsection 1.2 and
    Section 3) and trace identities (subsection 1.3 and Section 4).
    The  two pairs of rank 2 (see (1.4) and (1.5)) are chosen  so that many known identities
     of 3 terms given by Long [12] can be
     recovered
       (subsections 3.3 and 4.4).
       While  matrix identities (M1)-(M5)  are  known,
        trace identities (T1)-(T7)
      seem to be less circulated.
                 To illustrate,
 we give  a trace  identity as follows ((T6) of  Table 5).
         $$
   L_{2m}^2+5F_{2n+1}^2 +5F_{2m+2n+1}^2= 5L_{2m}F_{2n+1} F_{2m+2n+1} +4     .\eqno(1.2) $$

\noindent
 See Section 5 for pairs of rank 3 and their            matrix identities.
In the case $z$ is singular, the  pair $S(z, v)$ such that $zv=vz$ gives interesting
 identities as well. (1.3) of the following
 is such an example (see (v) of Section 7).
$$F_nL_{m+1}+F_{n-1}L_{m-1}+F_{n+1}L_{m}= L_{n+m+1}.\eqno(1.3)$$

\subsection{ Fibonacci pairs of rank $s$}
 The first pairs of rank two  in this article is
    $a=F_r$ and
    $$ w = F_r^{-1}\left (\begin{array}{cc}
F_{r+1}&F_1\\
(-1)^{r+1} F_1&(-1)^{r+1} F_{1-r}\\
\end{array}\right),\,\,
x = F_r^{-1}\left (\begin{array}{cc}
L_{r}&L_{0}\\
(-1)^{r+1}L_{0}& (-1)^{r+1}L_{-r}\\
\end{array}\right).\eqno(1.4)$$

\noindent Note that  the trace of $x$ is 0.
  Note also  that the Fibonacci pair in Theorem 1.1 is a special case of (1.4).
  The second  pair  is
 $$a=2,\,\,w = 2^{-1}\left (\begin{array}{cc}
1& 5\\
1& 1\\
\end{array}\right),\,\,
x = 2^{-1}\left (\begin{array}{cc}
0& 10\\
 2& 0\\
\end{array}\right).\eqno(1.5)$$

\noindent The detailed construction of the above  pairs of rank 2 can be found
  in Section 3.
  See Sections 5 and 6 for the construction of
 $F_s(w,x)$, where $s\ge 3$.

\subsection {Matrix identities of Fibonacci pairs of rank $s$}
 Let $F_s(w, x)$ be a Fibonacci pair of rank $s$ and let
$M, N\in \Delta = \{ w^a, xw^b\}$.
  The product $MN$  gives an equation. Take the pair $(M,N) = (w^m, xw^n)$
   for instance, the product $MN$ gives the following equation.
    $$(w^m)(xw^n) = xw^{m+n}.\eqno(1.6)$$

 \smallskip
\noindent
 We call equation (1.6) the {\em matrix equation} of
 $(w^m, xw^n)$.
Such equation actually gives  identities
 of Fibonacci and Lucas numbers.
For instance, if $(w^n,xw^m)$ is given as in Theorem 1.1,
 then  equation  $w^n(xw^m)
     =xw^{m+n}$ gives
  the identity
  $$ F_{n+1}L_{m+1}+ F_nL_m= L_{n+m+1}.\eqno(1.7)$$

 \smallskip
 \noindent  To be more precise,    (1,1) entry of the
   following equation  after the product $(w^n)(xw^m)$ is
   simplified gives the identity
   $ F_{n+1}L_{m+1}+ F_nL_m= L_{n+m+1}.$
 $$ w^n(xw^m) =\left (\begin{array}{cc}
F_{n+1}&F_n\\
F_n&F_{n-1}\\
\end{array}\right)
\left (\begin{array}{cc}
L_{m+1}&L_m\\
L_m&L_{m-1}\\
\end{array}\right)=
 \left (\begin{array}{cc}
L_{m+n+1}&L_{m+n}\\
L_{m+n}&L_{m+n-1}\\
\end{array}\right)=xw^{m+n}.\eqno(1.8)$$

\smallskip
  \noindent

\noindent Identity (1.7) is called a  {\em matrix identity}
 of  $(w^n,xw^m)$.
 Since $\Delta$ consists of 2 members, there are
  3 matrix equations  given as follows.

\smallskip
\begin{center}
Table  1  :  Matrix equations  of  $F_s(w,x)$
\end{center}
{ $$\begin{array}{llr}
  (M, N) && \mbox{Equations  associated with  } (M, N)\\
\\
(w^{n}, w^{m} )&&
 w^nw^m=w^{m+n}\\


\\
(w^{n}, xw^{m} )&&
w^{n}( xw^{m}) =xw^{n+m}
 \\


\\
(xw^{n}, xw^{m} )&&
(xw^{n})( xw^{m} )=x^2w^{n+m}
 \\


\end{array}
 $$}

\smallskip
\noindent
Matrix identities (such as (1.7)) coming from Table 1
are given in Sections 3 and 5. They are listed as
(M1)-(M5) and (N1)-(N5).

 \subsection{Trace identities of Fibonacci pairs of rank 2}
 Let $M$ and $N$ be given as in subsection 1.2.
   Set
   $A= M/\sqrt{d(M)}, \,\,B= N/\sqrt{d(N)}$,
   where $d(X)$ is the determinant of $X$.
    In the case $M$ and $N$ are $2\times 2$ matrix,
   since $A$ and $B$ commute with each other,  we can
    show in Section 4 that $ t(AB)+t(A^{-1}B)=t(A)t(B) $,
$t(BAB) +t(A)=t(B)t(AB)$ and
    $$t(A)^2+ t(B)^2 +t(AB)^2= t(A)t(B)t(AB) + 4, \eqno(1.9)$$

\noindent  where $t(X)$ is the trace of $X$.  These 3 identities
 are called the {\em trace identities}
  of $(M,N)$. See Section 4 for trace identities of $F_2(w, x)$.
  They are listed as (T1)-(T11).
   The following
        is a combination of a trace identity ((T1) of Table 5) and
   a matrix identity ((M4) of Section 3).
       $$L_{2n}^2 + L_{2m}^2+5F_{2m+2n}^2
 =   L_{2n}L_{2m}L_{2m+2n}.\eqno(1.10)$$

\subsection{Discussion} Our  goal is to relate
   identities  to
   matrices.
    We are far  from being done. For instance, we are unable
      to construct the identity
      $F_{3n}= F_{n+1}^3+F_n^3-F_{n-1}^3$ directly by
       Fibonacci pairs.
        The closest we can get is ((N5) of Table 7)
   $$ 3F_{n-1}^3-F_n^3+F_{n+1}^3+F_{n+1}F_nL_{n-1}
=2L_nF_{2n-1}.\eqno(1.11)$$

\smallskip
\noindent {\em An advantage of our method is that one does not have to verify
 the truth of the  identities  as they come automatically
  from the multiplication of  matrix equations}
   (see (1.8)). {\em Another advantage is that
    Fibonacci pairs can be constructed easily} (Propositions 2.1 and 5.1).

\subsection{Notations and organisation} Section 2
 explains how  $w$ and $x$ are constructed for $F_2(w, x)$.
    Section 3 gives two Fibonacci  pairs of rank 2 and
    their  matrix identities. It is shown that
    identities (13)-(24) of Long [12]
     are matrix  identities (subsection 3.3).
      Section 4 studies trace identities of $F_2(w, x)$. It is proved that identities (5)-(9) and (11)
     of  Hoggatt, Jr. and Bergum [7]  are  trace identities (subsection 4.4).
       Sections 5 and 6 study $F_s(w, x)$, where $s\ge 3$.
       $F_n$ is the $n$-th Fibonacci number and $L_m$ is the $m$-th Lucas numbers. Recall that
     $$F_0=0, F_1=1, F_{-r}= (-1)^{r+1}F_r,\,\,\,
       L_0=2, L_1=1, L_{-r}=(-1)^rL_r.
                       \eqno(1.12)$$

 \noindent The equations in Table 1 are called the
  {\em matrix equations} of  $F_s(w,x)$ and
   the identities (such as (1.7)) are called the
    {\em matrix identities} of
     $F_s(w, x)$.


\subsection{ Historical background} Matrices have been used to
 construct Fibonacci and Lucas identities. See [2]-[5], [9]
  [13],
  [14] and [15] for some detailed investigation.
  A  pair of matrices similar to (3.4) has been studied by
   Prasanta Kumar Ray [13]. To the best of our knowledge,
our systematic study of   Fibonacci pair $F_e(w, x)$
 (such that $\left < \bar w, \bar
x\right > \subseteq GL(s, \Bbb C)/\left < cI_s\,:\, c\ne 0\right > $  is not cyclic) is not in the literature yet.

\section{Fibonacci pairs of rank 2}
 We give detailed construction of Fibonacci  pairs of rank 2 in Proposition 2.1.  See
 Sections 5 and 6 for Fibonacci pairs of ranks 3 or more.

\subsection{} Let
 $F_2(w, x)$ be given as in subsection 1.1. Direct calculation  shows that
 the  characteristic polynomials
   of $w$ and $x$ are  $X^2-X-1$ and $X^2-5$ respectively. As a consequence,
$$w^n=F_nw+F_{n-1}I_2.\eqno(2,1)$$

  \subsection{}

   The existence of a Fibonacci pair of rank 2 is assured  by the
    following proposition.

   \smallskip
     \noindent {\bf Proposition 2.1.}  {\em Let $w$ be
      a  rational matrix whose characteristic polynomial
       is $X^2-X-1$ and let $a, b$ be rational numbers. Then
      $w^n= F_nw+F_{n-1}I_2 \mbox { and }
 wx=xw \mbox{ where }  x= aw+bI_2.$}

\smallskip \noindent
Proposition 2.1 can be proved by induction. Further,
  the proposition indicates that  the construction of a Fibonacci
 pair $F_2(w,x)$ is easy.
    For instance, $F_2(w_d, x(a,b))$
     is a Fibonacci pair, where
   $w_d={\tiny\left (\begin{array}{cc}
1-d& 1\\
 1+d-d^2&d\\
\end{array}\right)}$, $x(a,b)= aw_d+bI_2.$  In general, one can always take $w$ to be
 a matrix similar to the rational canonical form of $X^2-X-1$ over $\Bbb Q$.
To finish the subsection, we give a simple  example as follows. Set  $w= w_{-2}$. By (2.1), one  has
$$
  w^n =\left (\begin{array}{cc}
 L_{n+1}& F_n\\
 -5F_n&-L_{n-1}\\
\end{array}\right),\,\,
 wx=xw
 \mbox{ where }
 x= x (1,2)= w+2I_2,
 \eqno(2.2)$$

\noindent
 Hence $F_2(w, x)$ is a Fibonacci pair of rank 2.
 Note that $(2.2) $  is different from $(1.4)$ and (1.5).
  Note also that (1.4) and (1.5) are special cases of Proposition 2.1.

\subsection{} Let $F_2(w, x)$ be given as in (1.4) and (1.5). The set of all
 matrices commute with $w$ is a vector space of dimension 2
  over $\Bbb C$.  $\{w, x\}$ is  a basis of $V$. Let  $F_2(w, y)$
   be another Fibonacci pair.
   Then $y= aw+bx$ for some $a, b \in \Bbb Q$.
     Consequently,  matrix
      equations of  $F_2(w, y)$ are just combinations
       of matrix equations of  $F_2(w, x)$. In conclusion, $F_2(w, x)$ gives all the matrix equations and identities.
        Note that $x$ has trace 0 and $x^2$ is a scalar matrix
               This  simplifies the calculation
          of the last  matrix equation of Table 1.
        Connections
           between trace identities (see (1.9)) of  $F_2(w, x)$
           and $F_2(w, y)$ are quite involved
            (see subsection 4.4).

\smallskip
Let $F_2(w,x)$ be a Fibonacci pair. See (vi) of Section 7 if $x^2$ is not a scalar matrix.

   \section{Matrix identities of Fibonacci pairs of rank 2}
      We shall   construct
       2 Fibonacci pairs of rank 2 and 5
        matrix identities (M1)-(M5). We will see in subsection 3.3 that
       identities (13)-(24) of Long [12] are consequences
        of  (M1)-(M5).

 \subsection{Matrix identities (M1)-(M3)}   Set $a= F_r$,
    $$ w = F_r^{-1}\left (\begin{array}{cc}
F_{r+1}&F_1\\
(-1)^{r+1}F_1&(-1)^{r+1} F_{1-r}\\
\end{array}\right),\,\,
x = F_r^{-1}\left (\begin{array}{cc}
L_{r}&L_{0}\\
(-1)^{r+1}L_{0}& (-1)^{r+1}L_{-r}\\
\end{array}\right).\eqno(3.1)$$

\noindent
 It is clear  that $wx=xw$,
   trace of $x$ is 0  and  $x^2$ is a scalar matrix.
  By (2.1) and the known identity
   $F_{n+r}= F_nF_{r+1}+F_{n-1}F_r$, one has
{\small $$ w^n = F_r^{-1}\left (\begin{array}{cc}
F_{n+r}&F_n\\
(-1)^{r+1}F_n&(-1)^{r+1} F_{n-r}\\
\end{array}\right),\,\,
xw^n = F_r^{-1}\left (\begin{array}{cc}
L_{n+r}&L_{n}\\
(-1)^{r+1} L_{n}& (-1)^{r+1}L_{n-r}\\
\end{array}\right).\eqno(3.2)$$}

\smallskip
\noindent Hence $F_2(w,x)$ in (3.1) is a Fibonacci pair.
Note that $w$ and $w^n$ are closely related to
 $F_n^2 -F_{n+r}F_{n-r} = (-1)^{n-r} F_r^2$.
 Applying the technique we presented
 in subsection 1.2 (see (1.7) of subsection 1.2), we obtain 3  matrix identities. {\em They are direct consequences of $(3.2)$ and
  matrix multiplications }(see subsection 1.4 and Table 1).
  We tabulate our results in the following table.

\smallskip
\begin{center}
Table 2 :  Matrix identities for $F_2(w,x)$ given as in (3.1)
\end{center}
\smallskip
$$
\begin{array}{llrrr}
 (M, N) && \mbox{Matrix identity associated with  } (M, N) && \mbox{Name}\\
\\
(w^{n}, w^{m} )&&
 F_{n+r}F_{m+r}+(-1)^{r+1} F_nF_m=F_r F_{n+m+r} && (M1) \\


\\
(w^{n}, xw^{m} )&&
 F_{n+r}L_{m+r}+(-1)^{r+1}F_n L_m= F_rL_{n+m+r} && (M2)\\


\\
(xw^{n}, xw^{m} )&&
 L_{n+r}L_{m+r}+(-1)^{r+1}L_n L_m= 5F_{r}F_{n+m+r} && (M3)\\


\\

\end{array}
 $$
The table is read as follows. Let $(M,N)$ be given
 as in the first column. Then one can get an identity
  given as in the second column
  by considering the multiplication of $M$ and $N$.
  Take the first row for instance, one has
  $ (w^{n})(w^{m}) = w^{n+m}$
     and  (1,1)
   entry of this equation gives
$F_{n+r}F_{m+r}+(-1)^{r+1}F_nF_m=F_r F_{n+m+r}$.

\smallskip

\noindent {\bf Discussion.}
(i) Recall that $L_{-a}= (-1)^aL_a$ and $F_{-a}= (-1)^{a+1}F_a$.

\smallskip
\noindent
  (iii)
 Since  $wx=xw$ and the characteristic polynomials of $w$ and $x$ are  $X^2-X-1$
   and $X^2-5$, the traces
  of $w^n$ and $xw^n$ are  $L_n$ and $5F_n$
   respectively. Since   traces of $w^n$ and $xw^n$ can be calculated by (3.2) also,  one has
  $$
  F_{n+r}+(-1)^{r+1}F_{n-r}= F_rL_n\,\mbox { and  }\,
  L_{n+r}+(-1)^{r+1}L_{n-r}= 5F_rF_n.\eqno(3.3)$$

\subsection{Matrix identities (M4) and (M5)} To finish our study
 of Section 3, we consider
 $$ a=2,\,\,w=2^{-1} \left (\begin{array}{cc}
1&5\\
1&1\\
\end{array}\right)
,\,\, x= 2^{-1}\left (\begin{array}{cc}
0&10\\
2&0\\
\end{array}\right)
.\eqno(3.4)$$

\smallskip
\noindent  It is clear  that
 $xw=wx$ and  $x^2$ is a scalar matrix.
  By (2.1), $w^n$ and $xw^n$ are given as follows.
$$ w^n = 2^{-1}\left (\begin{array}{cc}
L_n&5F_n\\
F_n& L_n\\
\end{array}\right),\,\,
xw^n= 2^{-1}\left (\begin{array}{cc}
5F_{n}&5L_{n}\\
L_{n}& 5F_{n}\\
\end{array}\right).\eqno(3.5)$$

\smallskip
\noindent
 Hence $F_2(w,x)$ in (3.4) is a Fibonacci pair.  Note that the idea of the
   construction of $w$  comes from the fact that
    $L_n^2-5F_n^2= 4(-1)^n.$
     Since  the
  first  matrix equation  and the last  matrix equation (see Table 1) give the
    same matrix identities when $F_2(w, x)$ is given as in (3.4),
     we list 2 rather than 3 identities.
  They are direct consequences of $(3.5)$ and
  matrix multiplications (see subsection 1.4 and Table 1).

\medskip
\begin{center}
Table 3 :  Matrix identities for $F_2(w,x)$ given as in (3.4)

\end{center}
 \smallskip
 $$
\begin{array}{llrrr}
 (M, N) && \mbox{Matrix identity associated with  } (M, N)&&\mbox{Name}\\
\\
(w^{n}, w^{m} )&&
 L_nL_m+5F_nF_m = 2L_{m+n} && (M4)\\

\\

(w^{n}, xw^{m} )&&
L_nF_m +F_nL_m=2F_{m+n} && (M5)
\\




\end{array}
 $$

    \subsection{Discussion} In [12], Long exhibited
     17 identities (labeled as (13)-(24)) by studying $L_n^2-5F_n^2
      =4(-1)^n$.
     The main theme of this subsection is to
         show that

 \smallskip      \noindent {\bf Proposition 3.1.} {\em
  Identities $(13)$-$(24)$ of Long $[12]$ are matrix identities
   or combination of matrix identities $(M1)$-$(M5)$.}

        \smallskip
        \noindent {\em Proof.} One  sees that 8 of the
         17 identities are listed in  $(M1)$-$(M5)$.
          They are (14 even), (15 even), (16 even), (18), (19), (21),
           (22), and (23).
          The remaining 9 identities are combinations of
                   identities
          $(M1)$-$(M5)$. Take (24) and (20) of Long [11] for examples.
           The identities are
               $$
          F_nL_m-L_{n-d}F_{m+d}=(-1)^{n+1}L_{-d}F_{m-n+d},\,
          L_mL_n-5F_{n-d}F_{m+d}=(-1)^nL_{-d}L_{m-n+d}
          .\eqno(3.6)$$

  \smallskip
  \noindent      The first identity is a combination of $
  (M2)$, and  $(M5)$ and the second  is a combination of
  $(M1)$ and $(M4)$. The other 7 identities can be
   obtained similarly.\qed

\section{Trace identities of Fibonacci pairs of rank 2 }
 We will prove three  well known  trace identities
  ((4.2), (4.5A) and (4.5B)) and apply
 such identities to get trace identities of
  Fibonacci and Lucas numbers.

 \subsection{Trace function}
 A matrix $X$ is unimodular if its determinant $d(X)$ equals 1.
 Let $A$ and $B$
 be $2\times 2$ unimodular  matrices. Then
 $$t(A)^2 + t(B)^2 +t(AB)^2= t(A)t(B)t(AB) +
 t(ABA^{-1}B^{-1})+2,\eqno(4.1)$$

 \noindent
  where $t(X)$ is the trace of $X$. Note that Identity (4.1)
    holds only if $A$ and $B$ are $2\times 2$ matrices.
 In the case  $A$ and $B$ commute with each other, one has

 \smallskip
 \noindent {\bf Proposition 4.1.} {\em Let $A$ and $B$ be
  $2\times 2$  unimodular matrices such that $AB=BA$.
  Then
 $$t(A)^2+ t(B)^2 +t(AB)^2= t(A)t(B)t(AB) + 4,\eqno(4.2)$$

 }

 \smallskip
 \noindent {\em Proof.}
  Since $AB=BA$,
    $A$ and $B$ can be simultaneously upper-triangulated.
     Since the trace of a matrix is invariant under
      the change of basis,
   we may assume that
   $$A= \left (\begin{array}{cc}
a&*\\
0&a^{-1}\\
\end{array}\right) \mbox { and that }
 B= \left (\begin{array}{cc}
b&*\\
0&b^{-1}\\
\end{array}\right).\eqno(4.3)$$

\smallskip
\noindent Identity (4.2) now follows by straightforward calculation.
 Note that $t(ABA^{-1}B^{-1})=2$.
  \qed

\smallskip

\subsection{Trace identities (T1)-(T7)}
 Let  $F_2(w,x)$
 be given as
  in  (3.1) or  (3.4)  and let
 $M, N \in \{w^a,xw^b\}$. We may  apply (4.2)
    to $A=M/\sqrt {d(M)}$ and
   $ B= N /\sqrt {d(N)}$
    to get trace identities.
     As such identities have great resemblance,
      we  present
     7 such identities only in the following two tables.

\medskip

\begin{center}
Table 4 :  Trace identity (4.2) for $F_2(w,x) $ given as in (3.1)
\end{center}
\smallskip
{
$$
\begin{array}{lrr}
 (M, N) & \mbox{Trace identity (4.2) associated with  } (M, N)&\mbox{Name}\\

  \\
  (w^{2n}, w^{2m} )&
  L_{2n}^2 + L_{2m}^2+L_{2m+2n}^2
   =   L_{2n}L_{2m}L_{2m+2n} +4 &(T1)
  \\

\\
(xw^{2n}, xw^{2m })&
-5F_{2n}^2 -5F_{2m}^2+L_{2m+2n}^2
 = 5  F_{2n}F_{2m}L_{2m+2n} +4 &(T2)
\\



\\
(xw^{2n}, xw^{2m +1})&
-5F_{2n}^2 +5F_{2m+1}^2-L_{2m+2n+1}^2
 = -5  F_{2n}F_{2m+1}L_{2m+2n+1} +4 &(T3)
\\

\\
(xw^{2n+1}, xw^{2m+1} )\,\,&
5F_{2n+1}^2 +5F_{2m+1}^2+L_{2m+2n+2}^2
 = 5  F_{2n+1}F_{2m+1}L_{2m+2n+2} +4 &(T4)
\\

\end{array}
 $$}
\noindent

 Note that $(T1)$ is a special case of the trace identity for
  $(w^n, w^m)$ given as follows.
$$
  (-1)^nL_n^2 +(-1)^m L_m^2+(-1)^{m+n}L_{m+n}^2
   = (-1)^{m+n}  L_nL_mL_{m+n} +4.\eqno(4.4)$$


\noindent Note also that
trace identities for $(w^n, xw^m)$ and  $(xw^n, xw^m)$
 can be obtained also. However, we  do not include them here as such identities take more space to present.

\medskip

\begin{center}
Table 5 :  Trace identity (4.2)  for $F_2(w,x) $ given as in (3.4)
\end{center}
\smallskip
$$
\begin{array}{lrr}
 (M, N) & \mbox{Trace identity (4.2) associated with  } (M, N)&\mbox{Name}\\

\\
(w^{2n}, xw^{2m} )&
L_{2n}^2 -5F_{2m}^2-5F_{2n+2m}^2=-5L_{2n}F_{2m}F_{2n+2m}+4
 &(T5)
\\

\\
(w^{2n}, xw^{2m+1} )&
L_{2n}^2+ 5F_{2m+1}^2+5F_{2n+2m+1}^2=5L_{2n}F_{2m+1}F_{2n+2m+1}+4
 &(T6)
\\

\\
(w^{2n+1}, xw^{2m} )\,\,&
-L_{2n+1}^2- 5F_{2m}^2+5F_{2n+2m+1}^2=5L_{2n+1}F_{2m}F_{2n+2m+1}+4
 &(T7)
\\

\end{array}
 $$
\smallskip

\subsection{Trace identities (T8)-(T11)} Let $A$ and $B$
   be two by two unimodular  matrices. Then
 $$ (t(AB)+t(A^{-1}B)=t(A)t(B)  \,\,  \,\,\,(4.5A)\,\,\, \mbox{ and }\,\,\,\,
t(BAB) +t(A)=t(B)t(AB).\eqno(4.5B)$$

\smallskip
\noindent Note that  (4.5A) and (4.5B) can be proved easily if $AB=BA$
 (see the proof of Proposition 4.1).
 Let $(M, N)$ be given as in subsection 4.2. We may apply (4.5A) and (4.5B)
    to $A=M/\sqrt {d(M)}$ and
   $ B= N /\sqrt {d(N)}$
    to get identities.
     As such identities have great resemblance,
      we present 4 such identities only
       in the following table.

\bigskip

\begin{center}
Table 6 :  Trace identities (4.5A) and  (4.5B)  for $F_2(w,x) $ given as in (1.1)
\end{center}

{
$$\begin{array}{lrrrrr}
(M, N) &\mbox{Identity (4.5B) for  } (M, N)&
 \mbox{Name} & \mbox{Identity (4.5A) for  }
 (M, N) &\mbox{Name}  \\

\\
( w^{n},w^{m} )&
L_{2m+n  } +(-1)^mL_n=L_mL_{m+n}& (T8)&
L_{m+n} +(-1)^nL_{m-n} =L_mL_n &(T9)
\\


\\
(x w^{n},w^{m} )&
F_{2m+n} +(-1)^mF_n=L_mF_{m+n}&(T10)&F_{m+n} -(-1)^nF_{m-n} =L_mF_n
 &(T11)
\\

\end{array}
 $$}

\subsection{Discussion.} (i)
 Identities (5)-(9)  and (11) of  Hoggatt, Jr. and Bergum [7]
 are  trace identities (see Table 6).
 (ii) Let $F_2(w,x)$ be given as
 in (3.1) and let $y = w+x$. Identity (4.2) for
 $(M, N) =(w^n, yw^m)$ where $m$ and $n$ are even takes the following form
 $$
-11 L_n^2+
 (L_{m+1}+5F_m)^2 +(L_{m+n+1}+5F_{m+n})^2
  = L_n
(L_{m+1}+5F_m)(L_{m+n+1}+5F_{m+n})-44.\eqno(4.6)$$
Identity (4.6) is a combination of
 $11L_n^2 -11L_{m+1}^2-11L_{m+n+1}^2=
 -11L_nL_{m+1}L_{m+n+1}+44$ (see (4.4))
  and a simple but slightly lengthy  identity.
  This is true in general.
 We leave it to the readers to determine whether one
  should  claim that  trace identities of $F_2(w, y)$ can be obtained from   trace
    identities of $F_2(w, x)$.
 We do not work on  trace identities associated
  with  $F_2(w, y)$ $(y = aw+bx)$  as such identities
   are lengthy (see (4.6) for instance)
     before
  simplification.

\section {Fibonacci pairs of rank 3}

We give two Fibonacci pairs $F_3(w,x)$ and $F_3(z, v)$ of rank 3.
 The choice of $w$ is closely related to the Pascal's
triangle and $z$ is the rational canonical form of a polynomial.

\subsection{Fibonacci pairs of rank 3}
          The existence of
            $F_3(w, x)$ is guaranteed by the following  proposition which can be proved by induction.

     \smallskip
     \noindent {\bf Proposition 5.1.}  {\em Let $w$ be
      a $3\times 3$ rational matrix whose characteristic polynomial
       is $f(x)= X^3-2X^2-2X+1$ and let $a, b,c$ be rational numbers. Then
      $$w^n= F_nF_{n-1}w^2+F_nF_{n-2}w - F_{n-1}F_{n-2}I_3 \mbox { and }
 wx=xw \mbox{ where }  x= aw^2+bw +cI_3.
 \eqno(5.1)$$
   Note that $f(x) $ is the
   second auxiliary polynomial of the Fibonacci
  sequence  $($see $[2]$, $[12])$.
   }

\smallskip
  \noindent
{\bf Example 5.2.}  Let $w, x, z, $ and $v$ be three by three matrices given as follows.
$${ w=
 \left (\begin{array}{ccc}
0&0&1\\
0&1&2\\
1&1&1\\
\end{array}\right),\,
x=
 \left (\begin{array}{ccc}
0&1&0\\
2 &1&2\\
0&1&2\\
\end{array}\right),\,
z=
 \left (\begin{array}{ccc}
0&0&-1\\
1&0&2\\
0&1&2\\
\end{array}\right),\,
v=
 \left (\begin{array}{ccc}
-2&0&-1\\
1&-2&2\\
0&1&0\\
\end{array}\right)}.\eqno(5.2)$$

\noindent Note that   $w$ and $z$ have
   the same characteristic polynomial   $X^3-2X^2-2X+1$,
    the characteristic polynomial of $x$
     is $(X-1)(X^2-2X-4)$ and $z$ is the rational canonical form
      of $X^3-2X^2-2X+1$. Note also that $wx=xw$ and that  $zv=vz$.
     Applying  (5.1), $w^n$  can be calculated easily.
 $$
{w^n =
 \left (\begin{array}{ccc}
\phantom{ \Big |}F_{n-1}^2& F_{n-1}F_n& F_n^2\\
\phantom{ \Big |}2F_{n-1}F_n&F_{n-1}^2+F_nF_{n+1}&2F_nF_{n+1}\\
\phantom{ \Big |}F_n^2&F_nF_{n+1}&F_{n+1}^2\\
\end{array}\right)},\,\,
t(w^n)=
2F_{n-1}^2+F_nF_{n+1}+F_{n+1}^2.
\eqno(5.3)$$

\noindent
Recall that $t(X)$ is the trace of $X$.
The matrix $xw^n$ can be calculated easily as well (see (5.4)).
 As a consequence,  $F_3(w, x)$ is a Fibonacci pair of rank 3.
$$xw^n={
\left (\begin{array}{ccc}
\phantom{ \Big |} 2 F_{n-1}F_{n}&F_{n-1}^2+F_nF_{n+1}&2F_nF_{n+1}\\
\phantom{ \Big |} F_{n-1}^2+F_n^2+F_{n+1}^2&
F_{n-1}^2+2F_{n-1}F_n+3F_nF_{n+1}
&F_n^2
+F_{n+1}^2+F_{n+2}^2\\
\phantom{ \Big |}2 F_nF_{n+1}&F_{n+2}^2-F_nF_{n+1}&2F_{n+1}F_{n+2}\\
\end{array}\right)}.\eqno(5.4)$$

\noindent
 It is easy to see that  trace of $xw^n$ is
  $ F_{n-1}^2+F_{n+1}F_{n+3}
 +4F_{2n}$.
It is  also easy to see that   $F_3(z, v)$ is  a Fibonacci pair
 of rank 3 as well (see (5.6) for $z^n$).
  \subsection {Matrix identities}
Equations such as $(w^n)(w^m)$
$= w^{n+m},$
$(w^n)(xw^m)=xw^{n+m}$ and $ (xw^n)(xw^m) = x^2w^{n+m}$
  give matrix identities.
    The following table gives 5 of them.


\smallskip
\begin{center}
Table 7 :  Matrix identity   for $F_3(w,x) $ given as in (5.2)

\end{center}
\smallskip
$${
\begin{array}{lcrr}
 (M, N) & \mbox{entry } (r,s)&\mbox{Matrix identity associated with entry } (r,s)
 &\mbox{Name}\\

\\
( w,xw^{n} )&  (1,2)& F_{n+2}^2-F_n^2=F_{2n+2} &  (N1)   \\

\\

( w,xw^{n} )&  (2,1)& F_{n+2}^2-F_{n-1}^2=4F_nF_{n+1} &  (N2)\\

\\

( w,xw^{n} )&  (2,2)& F_{n-1}^2+F_{n+2}^2=F_nL_{n-1}+F_{n+2}L_n &  (N3)\\

\\

( w,xw^{n} )&  (3,2)& F_{n+3}^2-F_{n+2}^2-2F_{n-1}^2=F_{2n+2}+2F_{2n} &  (N4)\\

\\

( w^n,xw^{n} )&  (1,1)& 3F_{n-1}^3-F_n^3+F_{n+1}^3+F_{n+1}F_nL_{n-1}
=2L_nF_{2n-1} &  (N5)\\

\end{array}}
 $$

 \smallskip
 \noindent
The first row is read as follows.  (1,2) entry of
 the matrix equation $w(xw^n)= xw^{n+1}$
 gives the identity $F_{n+2}^2 -F_n^2 = F_{2n+2}$.
  The remaining  4 rows can be interpreted  similarly.
   Combination of (N1) and (N4) gives
    $F_{n+3}^2-2F_{n+2}^2-2F_{n+1}^2+F_n^2=0$
     (Hoggatt,  Jr. and  Bicknell [6]).

\subsection{Trace identities} We shall first prove the following lemma
 about  traces.

\smallskip
\noindent {\bf Lemma 5.3.} {\em
Suppose that the characteristic polynomial of $A$ is  $X^3-2X^2-2X+1$
 and that $A^n= (a_{ij}).$ Then
 $a_{11}+a_{22}+a_{33} = L_{2n}+(-1)^n$.  Note that
$L_{2n}+(-1)^n = F_{3n}/F_n.$}

 \smallskip
 \noindent {\em Proof.}
 Eigenvalues  of $A$ are $-1$,
 $\tau= (3+\sqrt 5)/2$ and $\sigma=3-\sqrt 5)/2$.
  Consequently,  eigenvalues of $A^n$ are $(-1)^n$, $\tau^n$ and $\sigma^n$.
  Since $\tau^n+\sigma^n= L_{2n}$, trace of
    $A^n=(a_{ij})$ is $L_{2n}+(-1)^n$.
   Hence
    $a_{11}+a_{22}+a_{33} = L_{2n}+(-1)^n.$\qed

 \smallskip
\noindent
 {\bf Example 5.4.} We give two trace identities, one for $w^n$ and one for
  $z^n$.
 Characteristic polynomial of $w$ is $X^3-2X^2-2X+1$.
 By  Lemma 5.3 and (5.3), one has
    $$
     L_{2n}+ (-1)^n
    = 2F_{n-1}^2+F_nF_{n+1}+F_{n+1}^2= 2F_{n-1}^2 + F_{n+1}F_{n+2}.
   \eqno(5.5)$$

   \noindent
   Since  the characteristic polynomial of $z$ is $X^3-2X^2-2X+1$,
     By (5.1), $z^n$ is given as follows.
   $$z^n={
 \left (\begin{array}{ccc}
\phantom{ \Big |} -F_{n-1}F_{n-2}& -F_nF_{n-1}& -F_{n+1}F_n\\
\phantom{ \Big |}F_{n}F_{n-2}& F_{n+1}F_{n-1}& F_{n+2}F_n\\
\phantom{ \Big |}F_{n}F_{n-1}& F_{n+1}F_{n}& F_{n+2}F_{n+1}\\
\end{array}\right),}\,\,\,\,
{\small t(z^n)= 2F_{n+1}^2-F_{n-1}F_{n-2}.}
\eqno(5.6)$$

 \noindent
 Applying  Lemma 5.3 and (5.6),  we have the identity
     $L_{2n}+ (-1)^n=
 2F_{n+1}^2- F_{n-1}F_{n-2}.$

\smallskip


\noindent {\bf Example 5.5.}
  Eigenvalues
 of $x$ are 1, $1+\sqrt 5$ and $1-\sqrt 5$. Since $wx=xw$, one can
 show
  that
  eigenvalues of $xw^n$ are
   $(-1)^n, $ $(1+\sqrt 5)\tau^n$ and $(1-\sqrt 5)\sigma^n$.
    It follows that the  trace of $xw^n$
     is $(-1)^n +L_{2n} +5F_{2n}$.
     Similar to Example 5.4,  expression of $xw^n$
      (see (5.4)) gives
             $$F_{n-1}^2+
             F_{n+1}F_{n+3}  = (-1)^n+L_{2n}+F_{2n}.\eqno(5.7)$$

\noindent
Note that if $B= (b_{ij})$ is a conjugate of $xw^n$, then
$b_{11}+b_{22}+b_{33} = (-1)^n+L_{2n}+5F_{2n}.$


\section{Fibonacci pairs of rank 4 or more}
Let $a$ and $b$ be  roots of $X^2-X-1=0$
and $\Phi_r(X)= \prod _{i=0}^r(X-a^ib^{r-i})$
 (see [2], [12]).
 Let $w$ be
    the rational canonical
    form of $\Phi_r(X)$.
     Similar to Proposition 5.1, one can show that
    $F_{r+1}(w, x)$ is a Fibonacci pair of rank $r+1$,
     where $x= aI_{r+1}+\sum a_iw^i$ for some $a\ne 0$ (see [10]).
        We do not pursue our study of  pairs of  rank
      4 or more
      in this article
      as the identities will be lengthy.

 \section{Final comment}
 (i) We have demonstrated
   how identities can be constructed by studying Fibonacci pairs.
    Such   work can be viewed as an extension of [11]
     that enables us to generate identities systematically.
   We hope the readers can construct more identities
   by applying similar methods.

 \smallskip \noindent (ii) It is fairly easy to get new matrix identities
   other than (M1)-(M5).
    For instance,
             if $F_2(w, x)$  is given as in subsection 2.2, where $w=w_{-2}$,  $x=x(1,2)$, then (1,1) entries of
            $(w^n)(w^m)=w^{n+m}$ and  $(w^n)(xw^m)=x w^{m+n}$ give the following
            identities.
           $$
           L_{n+1}L_{m+1}-5F_nF_m= L_{n+m+1}\, \mbox{  and }\,
           L_{n+1}F_{m+2}-F_nL_{m+1} = F_{m+n+2}
           .\eqno(7.1)$$

\smallskip \noindent
Note that identities in (7.1) are not included in (M1)-(M5) of Section 3.

\smallskip
\noindent  (iii) Let $z^n$ be given as in (5.6).
 (1,1) and (2,1) entries of  $(z^n)(z^n)=z^{2n}$ (after simplification) give the following identities.
$$L_{n-1}F_{2n-1} = 2F_n^3-F_{n-1}F_{n-2}^2\, \mbox{ and  }\,
L_nL_{n-1}= 2F_{n-1}F_{n-2} +F_{n+2}F_n.\eqno(7.2)$$

  \smallskip
  \noindent
  (iv) Let $w$ be given as
   in Proposition 2.1.  Consideration of  a  pair of
    $2\times 2$ matrices $L_2(w,v)$
     such that
   $vw^m = w^{-m} v$ and  that entries of $aw^m$ and $avw^m$
    are  Fibonacci or Lucas numbers $(a\ne 0)$ seems
     less fruitful as
      $vw = w^{-1} v$ implies that $v$ and
       $vw^n$ have trace 0 ((4.5A)).
     Can one find pairs of $2\times 2$ matrices other than $F_2(w, x)$ and $L_2(w, v)$ that
      give nice identities ?

\smallskip
\noindent (v)  Matrix identities of
a  pair of matrices $S(z, v)$ such that  $z$ is singular
 and $zv=vz$ are also of interest. For instance, if $z$ and $v$
  are given as follows.
   Then (1,1) entry of $(z^n)(vz^m)=vz^{n+m}$ gives
$F_nL_{m+1}+F_{n-1}L_{m-1}+F_{n+1}L_{m}= L_{n+m+1}$.
$$z=
{
  \left (\begin{array}{ccc}
1&0&1\\
1&-1&0\\
0&1&1\\
\end{array}\right)},\,\,\, v = 2I_3+z,\,\,
z^n=
{\small  \left (\begin{array}{ccc}
F_n& F_{n-1}&F_{n+1}\\
F_{n-2}&F_{n-3}&F_{n-1}\\
F_{n-1}&F_{n-2}&F_n\\
\end{array}\right)}
.\eqno(7.3)$$

\noindent  Matrix $z$ is taken from [8]. Note
 that the characteristic polynomial
  of $z$ is $X^3-X^2-X$.
  In general, we may take
  $z$ to be a matrix similar to the rational canonical form
   of $X^3-X^2-X$.

\smallskip
\noindent (vi)
 Let $F_s(w, x)$ be a Fibonacci pair of rank $s$.
  In the case $x^2$ is not a scalar matrix,
one  may  also consider matrix equation
 of the form $(x^aw^b)(x^cw^d)= x^{a+c}w^{b+d}$ to get
  matrix identities. We leave it to the readers if he or she finds
   it interesting.

\bigskip

\noindent MSC2010 : 11B39, 11B83.

\medskip

\end{document}